\documentclass{article}
\usepackage{amsthm,amsmath, amssymb, amsfonts}
\newtheorem{theorem}{Theorem}
\newtheorem{lemma}{Lemma}

\newtheorem{proposition}{Proposition}

\newcommand*{\R}{\mathbb{R}}
\newcommand*{\Z}{\mathbb{Z}}
\newcommand*{\Zp}{\mathbb{Z}_{\geq 0}}

\newcommand*{\pset}{\mathcal{P}}
\newcommand*{\icone}{\mathrm{int.cone}}
\newcommand*{\cone}{\mathrm{cone}}
\newcommand*{\lattice}{\mathrm{lattice}}

\newcommand*{\crank}{\mathrm{cr}}

\title{On the Carath\'eodory rank of polymatroid bases}
\author{Dion Gijswijt\footnote{CWI and Dep. of Mathematics, Leiden University. Email: dion.gijswijt@gmail.com.}\ and Guus Regts\footnote{CWI, Amsterdam. Email: regts@cwi.nl.}}

\begin{document}
\maketitle
\begin{abstract}
In this paper we prove that the Carath\'eodory rank of the set of bases of a (poly)matroid is upper bounded by the cardinality of the ground set. 
\\[.5cm]
\end{abstract}

\section{Introduction}
Let $H\subseteq \R^n$ be a finite set and denote by 
\begin{equation}
\icone(H):=\{\lambda_1x_1+\cdots+\lambda_kx_k\mid x_1,\ldots,x_k\in H, \lambda_1,\ldots,\lambda_k\in \Zp\}
\end{equation}
the integer cone generated by $H$. The \emph{Carath\'eodory rank} of $H$, denoted $\crank(H)$, is the least integer $t$ such that every element in $\icone(H)$ is the nonnegative integer combination of $t$ elements from $H$. 

The set $H$ is called a \emph{Hilbert base} if $\icone(H)=\cone(H)\cap \lattice(H)$, where $\cone(H)$ and $\lattice(H)$ are the convex cone and the lattice generated by $H$, respectively. 

Cook et al.\cite{CookFonluptSchrijver} showed that when $H$ is a Hilbert base generating a pointed cone, the bound $\crank(H)\leq 2n-1$ holds. This bound was improved to $2n-2$ by Seb\H o \cite{Sebo}. In the same paper, Seb\H o conjectured that $\crank(H)\leq n$ holds for any Hilbert base generating a pointed cone. A counterexample to this conjecture was found by Bruns et al.\cite{Brunsetal}.

Here we consider the case that $H$ is the set of incidence vectors of the bases of a matroid on $n$ elements. In his paper on testing membership in matroid polyhedra, Cunningham \cite{Cunningham} first asked for an upper bound on the number of different bases needed in a representation of a vector as a nonnegative integer sum of bases. It follows from Edmonds matroid partitioning theorem \cite{Edmonds} that the incidence vectors of matroid bases form a Hilbert base for the pointed cone they generate. Hence the upper bound of $2n-2$ applies. This bound was improved by de Pina and Soares \cite{dePina} to $n+r-1$, where $r$ is the rank of the matroid. Chaourar \cite{Chaourar} showed that an upper bound of $n$ holds for a certain minor closed class of matroids.

In this paper we show that the conjecture of Seb\H o holds for the bases of (poly)matroids. That is, the Carath\'eodory rank of the set of bases of a matroid is upper bounded by the cardinality of the ground set. More generally, we show that for an integer valued submodular function $f$, the Carath\'eodory rank of the set of bases of $f$ equals the maximum number of affinely independent bases of $f$.    

\section{Preliminaries}
In this section we introduce the basic notions concerning submodular functions. For background and more details, we refer the reader to \cite{Fujishige,Schrijver}.

Let $E$ be a finite set and denote its power set by $\pset(E)$. A function $f:\pset(E)\to \Z$ is called \emph{submodular} if $f(\emptyset)=0$ and for any $A,B\subseteq E$ the inequality $f(A)+f(B)\geq f(A\cup B)+f(A\cap B)$ holds. The set
\begin{equation}
EP_f:=\{x\in \R^E\mid x(U)\leq f(U)\text{ for all $U\subseteq E$}\}
\end{equation}
is called the \emph{extended polymatroid} associated to $f$, and 
\begin{equation}
B_f=\{x\in EP_f\mid x(E)=f(E)\}
\end{equation}
is called the \emph{base polytope} of $f$. Observe that $B_f$ is indeed a polytope, since for $x\in B_f$ and $e\in E$, the inequalities $f(E)-f(E-e)\leq x(e)\leq f(\{e\})$ hold, showing that $B_f$ is bounded.

A submodular function $f:\pset(E)\to \Z$ is the rank function of a matroid $M$ on $E$ if and only if $f$ is nonnegative, nondecreasing and $f(U)\leq |U|$ for every set $U\subseteq E$. In that case, $B_f$ is the convex hull of the incidence vectors of the bases of $M$.

Let $f:\pset(E)\to \Z$ be submodular. We will construct new submodular functions from $f$. The \emph{dual} of $f$, denoted $f^*$, is defined by 
\begin{eqnarray}
f^*(U):=f(E\setminus U)-f(E).
\end{eqnarray}
It is easy to check that $f^*$ is again submodular, that $(f^*)^*=f$ and that $B_{f^*}=-B_f$. For $a:E\to \Z$, the function $f+a$ given by $(f+a)(U):=f(U)+a(U)$ is submodular and $B_{f+a}=a+B_f$. The \emph{reduction of $f$ by $a$}, denoted $f|a$ is defined by
\begin{equation}
(f|a)(U):=\min_{T\subseteq U}(f(T)+a(U\setminus T)).
\end{equation}
It is not hard to check that $f|a$ is submodular and that $EP_{f|a}=\{x\in EP_f\mid x\leq a\}$. Hence we have that $B_{f|a}=\{x\in B_f\mid x\leq a\}$ when $B_f\cap\{x\mid x\leq a\}$ is nonempty. We will only need the following special case. Let $e_0\in E$ and $c\in \Z$ and define $a:E\to \Z$ by
\begin{equation}
a(e):=\begin{cases}c&\text{ if $e=e_0$,}\\f(\{e\})&\text{ if $e\neq e_0$.}\end{cases}
\end{equation}
Denote $f|(e_0,c):=f|a$. If $x_{e_0}\leq c$ for some $x\in B_f$, we obtain  
\begin{equation}
B_{f|(e_0,c)}=\{x\in B_f\mid x(e_0)\leq c\}.
\end{equation}

Our main tool is Edmonds' \cite{Edmonds} polymatroid intersection theorem which we state for the base polytope.
\begin{theorem}\label{edmonds}
Let $f,f':\pset(E)\to \Z$ be submodular. Then $B_f\cap B_{f'}$ is an integer polytope.
\end{theorem}

We will also use the following corollary (see \cite{Edmonds}).
\begin{theorem}\label{idp}
Let $f:\pset(E)\to \Z$ be submodular. Let $k$ be a positive integer and let $x\in (kB_f)\cap \Z^E$. Then there exist $x_1,\ldots,x_k\in B_f\cap \Z^E$ such that $x=x_1+\cdots+x_k$.
\end{theorem}

\begin{proof}
By the above constructions, the polytope $x-(k-1)B_f$ is the base polytope of the submodular function $f'=x+(k-1)f^*$. Consider the polytope $P:=B_f\cap B_{f'}$. It is nonempty, since $\frac{1}{k}x\in P$ and integer by Theorem \ref{edmonds}. Let $x_k\in P$ be an integer point. Then $x-x_k$ is an integer point in $(k-1)B_f$ and we can apply induction.
\end{proof}

Important in our proof will be the fact that faces of the base polytope of a submodular function are themselves base polytopes as the following proposition shows.
\begin{proposition}\label{faces}
Let $f:\pset(E)\to \Z$ be submodular and let $F\subseteq B_f$ be a face of dimension $|E|-t$. Then there exist a partition $E=E_1\cup\cdots\cup E_t$ and submodular functions $f_i:\pset(E_i)\to \Z$ such that $F=B_{f_1}\oplus\cdots\oplus B_{f_t}$. In particular, $F$ is the base polytope of a submodular function. 
\end{proposition}
A proof was given in \cite{Schrijver}, but for convenience of the reader, we will also give a proof here.
\begin{proof}
Let $\mathcal{T}\subseteq \pset(E)$ correspond to the tight constraints on $F$:
$$
\mathcal{T}=\{U\subseteq E\mid x(U)=f(U) \text{ for all $x\in F$}\}.
$$
It follows from the submodularity of $f$ that $\mathcal{T}$ is closed under taking unions and intersections.
Observe that the characteristic vectors $\{\chi^A\mid A\in \mathcal{T}\}$ span a $t$-dimensional space $V$. 
Let $\emptyset=A_0\subset A_1\subset\cdots\subset A_{t'}=E$ be a maximal chain of sets in $\mathcal{T}$. We claim that $t'=t$. Observe that the characteristic vectors $\chi^{A_1},\ldots, \chi^{A_{t'}}$ are linearly independent and span a $t'$-dimensional subspace $V'\subseteq V$. Hence $t'\leq t$.

To prove equality, suppose that there exists an $A\in \mathcal{T}$ such that $\chi^A\not\in V'$. Take such an $A$ that is inclusionwise maximal. Now let $i\geq 0$ be maximal, such that $A_i\subseteq A$. Then $A_i\subseteq A_{i+1}\cap A\subsetneq A_{i+1}$. Hence by maximality of the chain, $A_{i+1}\cap A=A_i$. By maximality of $A$, we have $\chi^{A\cup A_{i+1}}\in V'$ and hence, $\chi^A=\chi^{A\cap A_{i+1}}+\chi^{A\cup A_{i+1}}-\chi^{A_{i+1}}\in V'$, contradiction the choice of $A$. This shows that $t'=t$.

Define $E_i=A_i\setminus A_{i-1}$ for $i=1,\ldots, t$. Define $f_i:\pset(E_i)\to \Z$ by $f_i(U):=f(A_{i-1}\cup U)-f(A_{i-1})$ for all $U\subseteq E_i$. We will show that
\begin{equation}
F=B_{f_1}\oplus\cdots\oplus B_{f_t}.
\end{equation}
To see the inclusion `$\subseteq$', let $x=(x_1,\ldots,x_t)\in F$. Then $x(A_i)=f(A_i)$ holds for $i=0,\ldots,t$. Hence  for any $i=1,\ldots,t$ and any $U\subseteq E_i$ we have
\begin{equation}
x_i(U)=x(A_{i-1}\cup U)-x(A_{i-1})\leq f(A_{i-1}\cup U)-f(A_{i-1})=f_i(U),
\end{equation}
and equality holds for $U=E_i$.

To see the converse inclusion `$\supseteq$', let $x=(x_1,\ldots,x_t)\in B_{f_1}\oplus\cdots\oplus B_{f_t}$. Clearly
\begin{equation}
x(A_k)=\sum_{i=1}^k x_i(E_i)=\sum_{i=1}^k (f(A_i)-f(A_{i-1}))=f(A_k),
\end{equation}
in particular $x(E)=f(E)$. To complete the proof, we have to show that $x(U)\leq f(U)$ holds for all $U\subseteq E$. Suppose for contradiction that $x(U)>f(U)$ for some $U$. Choose such a $U$ inclusionwise minimal. Now take $k$ minimal such that $U\subseteq A_k$. Then we have
\begin{eqnarray}
x(U\cup A_{k-1})&=&x(A_{k-1})+x_k(E_k\cap U)\nonumber\\
&\leq& f(A_{k-1})+f_k(E_k\cap U)=f(U\cup A_{k-1}).
\end{eqnarray}
Since $x(A_{k-1}\cap U)\leq f(A_{k-1}\cap U)$ by minimality of $U$, we have
\begin{eqnarray}
x(U)&=&x(A_{k-1}\cup U)+x(A_{k-1}\cap U)-x(A_{k-1})\nonumber\\
&\leq &f(A_{k-1}\cup U)+f(A_{k-1}\cap U)-f(A_{k-1})\leq f(U).
\end{eqnarray}
This contradicts the choice of $U$.
\end{proof}

\section{The main theorem}
In this section we prove our main theorem. For $B_f\subseteq \R^E$, denote $\crank (B_f):=\crank (B_f\cap \Z^E)$. 
\begin{theorem}\label{main}
Let $f:\pset(E)\to \Z$ be a submodular function. Then $\crank (B_f)=\dim B_f+1$.
\end{theorem}

We will need the following lemma.
\begin{lemma}\label{directsum}
Let $B_{f_1}, \ldots, B_{f_t}$ be base polytopes. Then $\crank(B_{f_1}\oplus\cdots\oplus B_{f_t})\leq \crank(B_{f_1})+\cdots+\crank(B_{f_t})-(t-1)$.
\end{lemma}
\begin{proof}
It suffices to show the lemma in the case $t=2$.

Let $k$ be a positive integer and let $w=(w_1,w_2)$ be an integer vector in $k\cdot(B_{f_1}\oplus B_{f_2})$.  Let $w_1=\sum_{i=1}^r m_ix_i$ and $w_2=\sum_{i=1}^s n_iy_1$, where the $n_i,m_i$ are positive integers, the $x_i\in B_{f_1}$ and $y_i\in B_{f_2}$ integer vectors. Denote 
\begin{eqnarray}
\{0,m_1,m_1+m_2,\ldots,m_1+\cdots+m_r\}&\cup&\nonumber\\
\{0,n_1,n_1+n_2,\ldots,n_1+\cdots+n_s\}&=&\{l_0,l_1,\ldots,l_q\},
\end{eqnarray}
where $0=l_0<l_1<\cdots<l_q=k$. Since $m_1+\cdots+m_r=n_1+\cdots+n_s=k$, we have $q\leq r+s-1$. For any $i=1,\ldots,q$, there exist unique $j,j'$ such that $m_1+\cdots+m_{j-1}<l_i\leq m_1+\cdots+m_j$ and $n_1+\cdots+n_{j'-1}<l_i\leq n_1+\cdots+n_{j'}$. Denote 
$z_i:=(x_j,y_{j'})$. We now have the decomposition $w=\sum_{i=1}^q (l_i-l_{i-1})z_i$.
\end{proof}

We conclude this section with a proof of Theorem \ref{main}.
\begin{proof}[Proof of Theorem \ref{main}.]
The inequality $\crank (B_f)\geq \dim B_f+1$ is clear. We will prove the converse inequality by induction on $\dim B_f+|E|$, the case $|E|=1$ being clear. Let $E$ be a finite set, $|E|\geq 2$ and let $f:\pset(E)\to \Z$ be submodular.

Let $k$ be a positive integer and let $w\in kB_f\cap \Z^E$. We have to prove that $w$ is the positive integer combination of at most $\dim B_f+1$ integer points in $B_f$. We may assume that  
\begin{equation}
\dim B_f=|E|-1.
\end{equation}
Indeed, suppose that $\dim B_f=|E|-t$ for some $t\geq 2$. Then by Proposition \ref{faces}, there exist a partition $E=E_1\cup\cdots\cup E_t$ and submodular functions $f_i:\pset(E_i)\to \Z$ such that $B_f=B_{f_1}\oplus\cdots\oplus B_{f_t}$. By induction, $\crank (B_{f_i})=\dim B_{f_i}+1$ for every $i$. Hence by Lemma \ref{directsum}
\begin{eqnarray}
\crank (B_f)&\leq&\crank (B_{f_1})+\cdots+\crank (B_{f_t})-(t-1)\nonumber\\
&=&\dim B_{f_1}+\cdots+\dim B_{f_t}+1=\dim B_f+1.
\end{eqnarray}

Fix an element $e\in E$. Write $w(e)=kq+r$ where $r,q$ are integers and $0\leq r\leq k-1$. Let $f'=f|(e,q+1)$.
By Theorem \ref{idp}, we can find integer vectors $y_1,\ldots, y_k\in B_{f'}$ such that $w=y_1+\cdots+y_k$. We may assume that $y_i(e)=q+1$ for $i=1,\ldots,r$. Indeed, if $y_i(e)\leq q$ would hold for at least $k-r+1$ values of $i$, then we would arrive at the contradiction $w(e)\leq (k-r+1)q+(r-1)(q+1)\leq kq+r-1<w(e)$. 

Let $f'':=f|(e,q)$. Denote $w':=y_1+\cdots+y_r$. So we have decomposed $w$ into integer vectors
\begin{eqnarray}
w'&\in &rB_{f'}=B_{rf'}\nonumber\\
w-w'&=&y_{r+1}+\cdots+y_k\in (k-r)B_{f''}=B_{(k-r)f''}.
\end{eqnarray}
We may assume that $r\neq 0$, since otherwise $w\in kF$, where $F$ is the face $B_{f''}\cap \{x\mid x(e)=q, x(E)=f(E)\}$ of dimension $\dim F\leq |E|-2$ (since $|E|\geq 2$). Then by induction we could write $w$ as a nonnegative integer linear combination of at most $1+(\dim F)<\dim B_f+1$ integer vectors in $B_{f''}\subseteq B_f$. 

Consider the intersection
\begin{equation}
P:=B_{rf'}\cap B_{w+(k-r)(f'')^*}.
\end{equation}
Observe that $P$ is nonempty, since it contains $w'$. Furthermore, by Theorem \ref{edmonds}, $P$ is an integer polytope. Hence taking an integer vertex $x'$ of $P$ and denoting $x'':=w-x'$, we have that $x'$ is an integer vector of $B_{rf'}$ and $x''$ is an integer vector of $B_{(k-r)f''}$. 

Let $F'$ be the inclusionwise minimal face of $B_{rf'}$ containing $x'$ and let $F''$ be the inclusionwise minimal face of $B_{w+(k-r)(f'')^*}$ containing $x'$. Denote $H':=\mathrm{aff.hull}(F')$ and $H'':=\mathrm{aff.hull}(F'')$. Since $x'$ is a vertex of $P$, we have 
\begin{equation}
H'\cap H''=\{x'\}.
\end{equation}
Indeed, every supporting hyperplane of $B_{rf'}$ containing $x'$ also contains $F'$ by minimality of $F'$, and hence contains $H'$. Similarly, every supporting hyperplane of $B_{w+(k-r)(f'')*}$ containing $x'$ also contains $H''$. Since $x'$ is the intersection of supporting hyperplanes for the two polytopes, the claim follows. 

Observe that both $F'$ and $F''$ are contained in the affine space 
\begin{equation}
\{x\in\R^n\mid x(E)=rf(E),\ x(e)=r(q+1)\},
\end{equation}
which has dimension $n-2$ since $|E|\geq 2$. It follows that 
\begin{eqnarray}
\dim F'+\dim F''&=&\dim H'+\dim H''\nonumber\\
&=&\dim(\mathrm{aff.hull}(H'\cup H''))+\dim(H'\cap H'')\nonumber\\
&\leq& n-2.
\end{eqnarray}

Since $F''$ is a face of $B_{w+(k-r)(f'')^*}$ containing $x'$, we have that $w-F''$ is a face of $B_{(k-r)f''}$ containing $x''$. By induction we see that    
\begin{eqnarray}
\crank (F')+\crank (w-F'')&\leq& (\dim F'+1)+(\dim (w-F'')+1)\nonumber\\
&=&\dim F'+\dim F''+2\leq n.
\end{eqnarray}  
This gives a decomposition of $w=x'+x''$ using at most $n$ different bases of $B_f$, completing the proof.
\end{proof}


\end{document}